\documentclass[12pt]{article}

\usepackage{amssymb}
\usepackage{multicol}
\usepackage{turnstile}
\usepackage{lineno}
\modulolinenumbers[5]
\usepackage{caption}
\usepackage[a4paper]{geometry}
\usepackage{amsfonts}
\usepackage{arydshln}
\usepackage{xtab}
\usepackage{setspace}
\usepackage{fancyhdr}
\usepackage[dvips]{graphicx}
\usepackage[T1]{fontenc}
\usepackage[latin1,utf8]{inputenc}
\usepackage[main=english,italian]{babel}
\usepackage{type1ec}
\usepackage[final]{microtype}
\usepackage{amsmath,amsthm,eucal,dsfont,bm,mathrsfs,stmaryrd}
\usepackage{lmodern}
\usepackage{textcomp}
\usepackage{pict2e,enumitem}
\usepackage{hyperref}
\usepackage[dvipsnames]{xcolor}
\usepackage[nodayofweek]{datetime} 
\usepackage{comment}

\hypersetup{
    pdfpagemode={UseOutlines},
    bookmarksopen,
    pdfstartview={FitH},
    colorlinks,
    linkcolor={black},
    citecolor={black},
    urlcolor={black}
            }

\definecolor{cornellred}{rgb}{0.7, 0.11, 0.11}

\usepackage{geometry}

\newdateformat{monthyear}{\monthname[\THEMONTH], \THEYEAR}

\renewcommand{\P}{\mathbb{P}}

\theoremstyle{plain}
        \newtheorem{theorem}{Theorem}[]
        
        \newtheorem{corollary}{Corollary}[]
        
\theoremstyle{definition}

\renewcommand
        {\thefootnote}{\arabic{footnote}}
        
\newcommand{\symfootnote}[1]{%
\let\oldthefootnote=\thefootnote%
\stepcounter{mpfootnote}%
\addtocounter{footnote}{-1}%
\renewcommand{\thefootnote}{\fnsymbol{mpfootnote}}%
\footnote{#1}%
\let\thefootnote=\oldthefootnote%
}

\usepackage{tikz}

\makeatletter
\def\bbibitem#1{\item[]%
    \if@filesw\immediate\write\@auxout{\string \bibcite {#1}{\the\value{\@listctr }}}\fi\ignorespaces}
\makeatletter

\title{Confirming Mathematical Conjectures by Analogy
}
\author{Francesco Nappo\footnote{Politecnico di Milano, Department of Mathematics, via Bonardi 9, Campus Leonardo, 20133, Milan (Italy). E-mail: \texttt{francesco.nappo@polimi.it}}\ , Nicolò Cangiotti\footnote{Politecnico di Milano, Department of Mathematics, via Bonardi 9, Campus Leonardo, 20133, Milan (Italy). E-mail: \texttt{nicolo.cangiotti@polimi.it}} \ and Caterina Sisti\footnote{Università degli Studi di Torino, Department of Philosophy and Education Sciences, Palazzo Nuovo, Via Sant'Ottavio, 20, 10124 Turin (Italy). E-mail: \texttt{caterina.sisti@unito.it}} }
\date{}

\begin{document}

\maketitle

\begin{abstract}
Analogy has received attention as a form of inductive reasoning in the empirical sciences. However, its role in pure mathematics has received less consideration. This paper provides an account of how an analogy with a more familiar mathematical domain can contribute to the confirmation of a mathematical conjecture. By reference to case-studies, we propose a distinction between an incremental and a non-incremental form of confirmation by mathematical analogy. We offer an account of the former within the popular framework of Bayesian confirmation theory. As for the non-incremental notion, we defend its role in rationally informing the prior credences of mathematicians in those circumstances in which no new mathematical evidence is introduced. The resulting ‘hybrid’ framework captures many important aspects of the use of analogical inference in the realm of pure mathematics. 
\end{abstract}

\section{Introduction}
\label{Sec1}
Proof is the golden standard for mathematics. But proofs are also hard to find. Frequently, working mathematicians seek relief in evidence that merely supports a mathematical conjecture. Confirmation by instances is an example. According to Collatz’s (2010) conjecture, any process that consists in the application of a simple procedure (viz., to divide any input number by two if it is even or triple it and add one if the input is odd) will always end with the number $1$ (where the algorithm always ends). As of 2022, this yet unproven mathematical conjecture has been computer-checked for all starting integers up to $2^{68}$. This fact alone does not guarantee truth: several important numerical conjectures proved false for large enough values despite the enumerative evidence initially supporting them. In practice, however, mathematicians often regard the vastness of the sample of positive instances as non-negligible support for a conjecture.  

\medskip

This paper addresses another mode of non-deductive reasoning in mathematics: inference from analogy. This is the form of reasoning whereby, from the known similarities between the two mathematical domains, one concludes that some claim that holds in the more familiar source is also true of the target. In spite of its defeasible nature, analogical reasoning plays a prominent role in mathematical practice. For instance, one of the greatest results of 20th-century mathematics, due to Deligne (1973), proves a conjecture by André Weil that is a close analogue in algebraic geometry of Riemann’s hypothesis (or, more precisely, an analogue of the hypothesis for varieties over finite fields). According to several experts, Weil’s and Deligne’s results “provide some of the best reasons for believing that the Riemann hypothesis is true” (Edwards 1974:298) - even better than the enumerative evidence uncovered so far (cf. Deninger 1994:493).   

\medskip

Notwithstanding the central evidential role that analogies occupy in the practice of mathematicians, the mechanism whereby results gathered in a familiar mathematical source may support conjectures in an analogous target has not received nearly as much attention as the corresponding issue in the empirical sciences. We find that this is a significant omission. A central problem in the philosophy of mathematics is the question of how mathematical knowledge is possible (cf. Benacerraf 1965). A proper epistemology of mathematics should not neglect the actual research processes by which working mathematicians make new discoveries.

\medskip

The two most relevant discussions in the recent literature, by Corfield (2003:$\S$5.4) and Bartha (2009:$\S$5.7), make widely divergent suggestions on the matter. Drawing from Polya’s (1954b) pioneering work on analogy and induction in mathematics, Corfield (2003:103) defends Bayesian confirmation theory (BC) as a promising framework for developing an account of inductive support by mathematical analogy - without, however, providing full details as to how BC’s formal tools can be applied. Bartha (2009:279), on the contrary, denies that a Bayesian approach is even an option. On his view, an analogical argument in mathematics should not be seen as functioning as ordinary Bayesian evidence  - in the ‘incremental’ sense that the conjecture’s probability is greater in light of it than given the background knowledge alone; rather, it should be understood as influencing a mathematician’s prior credences before further evidence is sought. Bartha (2009:280) tentatively proposes to justify this role of analogy in establishing a ‘non-incremental’ form of confirmation for a conjecture (equivalent to regarding it as ‘\emph{prima facie plausible}’) by appealing to an alleged constraint of symmetry on rationality.

\medskip

Our aim in this paper is to make progress on the existing contributions, by outlining a novel systematic proposal concerning how analogies with more familiar domains may contribute to the confirmation of mathematical conjectures.\footnote{We are concerned here specifically with the descriptive issue of providing an account of how analogies actually underpin belief in mathematical conjectures, as opposed to the normative issue of whether this confirmatory use of analogy is indeed justified. We are skeptical that an answer to the latter problem will overturn the appearance that mathematicians do, as a matter of fact, rely on analogy to underwrite their confidence in conjectures (as Baker 2007 argues for the case of enumerative induction in number theory), but we will leave the matter to a separate discussion.} On our account - a middle ground between Corfield’s and Bartha’s - confirmation by analogy possesses both incremental and non-incremental aspects. The former, we will argue, has a relatively straightforward Bayesian representation as an instance of the phenomenon of transitivity of confirmation (cf. Hesse 1970; Roche and Shogenji 2013). As for the non-incremental aspect, we believe that a Bayesian treatment is not forthcoming. However, we will defend the role of the non-incremental notion in rationally informing the prior credences of mathematicians in those circumstances in which no new mathematical evidence is introduced. As we will discuss by reference to case studies from mathematical research, our framework captures important aspects of the logic of analogical inference in mathematics without invoking any suspicious symmetry constraint on rationality.

\medskip

The discussion below will proceed as follows. Section two will present the symmetry-based proposal advanced in Bartha (2009) and illustrate its systematic failures in recognizing the different forms of confirmatory support that analogical arguments can provide. Section three will move on to discuss the main obstacles for proposals, such as Corfield’s (2003), which purport to couch inductive support by mathematical analogy in fully Bayesian terms. In section four, we will offer detailed illustrations of how, assuming that some general issues with the application of BC to mathematics can be set aside, an analogical argument that introduces new mathematical evidence can confirm a conjecture from a Bayesian standpoint. In section five, we will address cases of inductive support by analogy that arguably escape the Bayesian definition of confirmation, while remaining compatible with the adoption of BC as a doctrine about rational credence updating in mathematics. Section six will conclude with a summary of our arguments.

\section{Not Just ‘Plausibility’}
\label{Sec2}

Bartha (2009) defends the view that an analogy with a more familiar domain can provide a form of inductive support to mathematical conjectures. However, unlike other notable defenders of this view (e.g., Polya 1954b), Bartha claims that inductive support from analogy should not be understood in the Bayesian terms of incremental confirmation - as when a piece of evidence rationally increases a hypothesis’ probability. His view is in two parts. First, Bartha claims “that plausibility of analogical reasoning is often best construed not precisely in terms of probabilities but in terms of \emph{relative conditional betting quotients}” (2009:182). Talk of betting quotients serves to escape the problem that our uncertainty towards conjectures is seemingly inexpressible in probabilistic terms because of \emph{logical omniscience}. The problem - which we will discuss in detail in section three - is that an unfettered application of the probability axioms apparently requires that one assigns a credence of either zero or one to every decidable mathematical claim.

\medskip

Second, Bartha claims that the “[conditional] betting quotients are \emph{symmetry based}” (2009:182). More specifically, Bartha’s view is as follows. Suppose that an agent assigns a relatively high conditional betting quotient to $H$, formally $Q (H \vert E)$, where $E$ is evidence about a mathematical source domain and $H$ is a well-confirmed hypothesis (or proven theorem) about that domain. She then finds out that $H$ has an analogue in another mathematical domain, $H^*$, that she had not considered before or such that her conditional betting quotient $Q(H^* \vert E^*)$ was negligibly low (where $E^*$ stands for a set of properties similar to those in E that the target is known to possess). Under those specific circumstances, Bartha (2009:182) contends that a symmetry constraint imposes a rational revision of her betting quotient for $H^*$ to some “non-negligible” value - viz., a value closer to $Q (H \vert E)$. The guiding idea is that there seems to be something asymmetric and incoherent about giving such a high quotient to $H$, without even giving a shot (so to speak) to the analogous bet in the target, absent decisive evidence against it.

\medskip

By means of the doctrine that probabilities about mathematical conjectures are betting quotients and that a general symmetry constraint governs them, Bartha recovers a sense in which an analogy with a familiar source can make a difference to a mathematician’s trust in a mathematical conjecture about a target without making use of probabilities. In particular, on Bartha’s account an analogy can justify assigning a ‘non-negligible’ value to one’s conditional betting quotients in some hypothesis H* before any new evidence for it is sought, equivalent to making the hypothesis “\emph{prima facie plausible}” (2009:296). Such a role in confirmation is eminently non-Bayesian in at least one sense: if Bartha is correct, we can have a rational revision of opinions regarding mathematical subjects in spite of the fact that the Bayesian definition of confirmation is not satisfied and that the revision is not mandated by the Bayesian rule of credence updating known as \emph{conditionalization} (whereby an agent’s posterior probability in a hypothesis must be exactly equal to its prior probability conditional on any new evidence).\footnote{As Bartha (2009:300) notes, the proposal escapes Dutch Book arguments insofar as the symmetry constraints is not an alternative rule of credence updating, but rather a synchronic constraint on credences in addition to coherence.} 

\medskip

In what follows, we will not question Bartha’s (2009) arguments for the claim that a symmetry principle constrains conditional betting quotients as described in his account. The concerns that we will voice below are exclusively directed at the applicability of Bartha’s framework to concrete examples of analogical reasoning in mathematics. That is to say, assuming for the sake of the argument that a symmetry-based justification can be made to work, is inductive support of mathematical conjectures from analogy best understood along the lines that Bartha indicates? We believe that the answer is negative. Let’s take a look at some of the examples that support our contention.

\subsection*{Example 2.1: Riemann’s Hypothesis}
\label{sec21}

One problem with Bartha’s symmetry-based account is that it has limited scope. It identifies a way for an analogy to contribute distinctly to the plausibility of a conjecture only when the latter is assigned no or very little conditional betting quotient to begin with. However, this neglects cases in which the conjecture that the analogy supports already enjoys ‘non-negligible’ credence (or betting quotient). For instance, it may be thought that Deligne’s (1973) proof of a close analogue of Riemann’s hypothesis in algebraic geometry adds to the credibility of the latter conjecture even though, before Deligne’s work, the hypothesis was already regarded as plausible by the majority of working mathematicians. Insofar as Bartha’s proposal only addresses those cases in which a mathematical conjecture was assigned no or negligible credence, it does not cover the role of a proof such as Deligne’s in indirectly supporting the Riemann hypothesis. 

\subsection*{Example 2.2: Euler Characteristic}
\label{sec22}

Second, the symmetry-based account seems unable to account for the different strengths that analogical arguments have in mathematical practice. The following example from solid geometry will be helpful to illustrate this point. Considering that in plane geometry the number of vertices ($V$) and edges ($E$) in convex polygons is the same ($V=E$), Euler (1758) asked if any similar regularity holds for the elements of all convex polyhedra. The correct answer is that a convex polyhedron’s edges ($E$), vertices ($V$) and faces ($F$) are governed by the following stable relation:
\begin{equation}
\label{e1}
V-E+F=2.
\end{equation}
Here are two distinct analogical arguments for this conclusion. The first (suggested by Polya 1954a:43) starts from the algebraic observation that the relation $V=E$ can be rewritten as:
\begin{equation}
\label{e2}
V-E+F=1.
\end{equation}
This suggests that the regularity we may be looking for is an alternating sum of the number of zero-dimensional elements ($V$), the number of one-dimensional elements ($E$), and that of two-dimensional elements ($F$). By analogy, and introducing the three-dimensional element of the number of solids ($S$), \eqref{e2} induces the following generalization to the three-dimensional case (which can be easily verified for cubes, tetrahedra, and dodecahedra):\footnote{Cf. Bartha (2009:155) for the claim that Polya’s reasoning is best construed in terms of an argument from analogy.}
\begin{equation}
\label{e3}
V-E+F-S=1.
\end{equation}
This reasoning yields the correct conclusion: since $S=1$ holds for all polyhedra, \eqref{e3} entails \eqref{e1}.

\medskip

An alternative argument, inspired to Cauchy’s (1811), runs as follows (cf. Lakatos 1976:7--9). Imagine to smash a polyhedron, a cube say, until it is flat. The edges do not break, but may be elongated in the process. The result is a two-dimensional figure with the same number of vertices and edges as the cube, but with one face less (one can think of it as two squares one inside the other with the corresponding vertices connected by an edge; see Fig. 1). If we now draw lines among the unconnected vertices of the smashed cube, we obtain small triangles of different sizes and shapes. We note that the relation between vertices, edges, and faces remains constant. Further, if we remove one of the triangles obtained by our drawing, there are only two options: either we remove an edge and a face, or we remove a vertex, a face and two edges. Either way, $V-E+F$ remains constant (i.e., it remains equal to $1$). By analogous reasoning, one conjectures that for any other ‘smashed’ polyhedron the relation $V-E+F$ will remain constant. Indeed, we can easily verify that the same is also true for tetrahedra and dodecahedra.\footnote{We think that analogy plays a central role in the reasoning in that the latter arguably derives its strength, not so much from the number of observed positive instances,(which are relatively few), but from the similarities that link all convex polyhedra together: specifically, from the fact that, because of the underlying commonalities that make them all instances of the same geometrical genus, all convex polyhedra appear to be amenable to the same smashing and triangulating operations to which we have subjected the cube, yielding exactly the same result. See also Section \ref{sec43}.}

\begin{figure}[ht]
\centering
\begin{minipage}[c]{.49\textwidth}
\centering
\begin{tikzpicture}[scale=1,>=latex]
  \draw[ultra thick] (0,0)  rectangle (5,5);
  \fill[fill=black,opacity=0.75] (3.5,3.5)--(5,5)--(5,0);
  \draw[ultra thick]  (1.5,1.5) rectangle (3.5,3.5);
  \draw[ultra thick] (0,0)--(1.5,1.5);
  \draw[ultra thick] (5,0)--(3.5,1.5);
  \draw[ultra thick] (0,5)--(1.5,3.5);
  \draw[ultra thick, dotted] (1.5,1.5) -- (3.5,3.5);
  \draw[ultra thick, dotted] (1.5,1.5) -- (0,5);
  \draw[ultra thick, dotted] (0,0) -- (3.5,1.5);
  \draw[ultra thick, dotted] (1.5,3.5) -- (5,5);
\end{tikzpicture}
\end{minipage}
\begin{minipage}[c]{.49\textwidth}
\centering
\begin{tikzpicture}[scale=1,>=latex]
  \draw[ultra thick] (0,0)  -- (5,0);
  \draw[ultra thick] (0,0)  -- (0,5);
  \fill[fill=black,opacity=0.75] (1.5,3.5)--(5,5)--(3.5,3.5);
  \draw[ultra thick]  (1.5,1.5) rectangle (3.5,3.5);
  \draw[ultra thick] (0,0)--(1.5,1.5);
  \draw[ultra thick] (5,0)--(3.5,1.5);
  \draw[ultra thick] (0,5)--(1.5,3.5);
  \draw[ultra thick, dotted] (1.5,1.5) -- (3.5,3.5);
  \draw[ultra thick, dotted] (1.5,1.5) -- (0,5);
  \draw[ultra thick, dotted] (0,0) -- (3.5,1.5);
\end{tikzpicture}
\end{minipage}
\caption*{Figure 1}
\end{figure}

The problem for Bartha’s account is that it treats the two arguments exactly alike in confirmatory potential, when they are not. The Polya-inspired reasoning is arguably of the weak heuristic variety, based on a rather superficial algebraic similarity between \eqref{e1} and \eqref{e2}. A mathematician who did not know that \eqref{e1} holds would be hardly moved by what appears to be a merely accidental algebraic similarity - one having no clear relation to the geometrical issue at stake.\footnote{We note that Polya (1954a) himself regarded the algebraic reasoning as a mere heuristic for generating conjectures.} Conversely, Cauchy’s geometrical reasoning is stronger in comparison.\footnote{Lakatos (1976:9) notes that several nineteenth-century mathematicians regarded Euler’s problem about the relation between vertices, edges, and faces in convex polyhedra to be practically settled by the Cauchy-inspired argument.} We cannot but feel that the imaginary operations involved in the reasoning are possible on any smashed polyhedron whatsoever. In each case, it is not clear what could go wrong to falsify \eqref{e1}. Yet, using Bartha’s approach no difference can be traced as to their respective capacity for inductive support: in both cases, we start with a high conditional betting quotient $Q (H \vert E)$, where $H$ is equation \eqref{e2} and $E$ includes all relevant facts about convex polygons; and in both cases, by symmetry, one ought to assign a ‘close enough’ value to the bet $Q(H^* \vert E^*)$, where $H^*$ is equation \eqref{e1} and $E^*$ includes all relevant facts about convex polyhedra (analogous to those in $E$). Hence, Bartha’s account clearly neglects the difference in strength between the two arguments.\footnote{ Bartha (2009) countenances the possibility that some analogical arguments are stronger than others as the result of the different credence assigned to the source theorem. However, in the chosen example we have the same source theorem for both analogical arguments, namely Equation \eqref{e2}. Hence, he cannot account for the difference in strength.}

\subsection*{Example 2.3: Area and Volume}
\label{sec23}

If the above example has not been sufficiently convincing, here is an even starker illustration that a symmetry-based account proves too much. It is well-known that for a rectangle of sides $x$ and $y$ the \emph{Area} can be computed by $x\cdot y$ and that, for a rectangular box of sides $x$, $y$ and $z$, the Volume of is given by $x\cdot y\cdot z$. Here is a strong analogical argument: given that, of all rectangles, the square maximizes \emph{Area}, it is plausible that, of all rectangular boxes, the cube maximizes \emph{Volume}. On Bartha’s view, the work of the analogical argument in making the conclusion plausible should be understood as follows. Let H be the theorem that the square maximizes \emph{Area} and $E$ the set of properties of two-dimensional rectangles that have analogues in the three-dimensional case. Since the conditional betting quotient $Q(H \vert E)$ is high, by symmetry a rational agent ought to assign a ‘non-negligible’ value to the conditional bet represented by $Q(H^* \vert E^*)$, where $H^*$ stands for the claim that the cube maximizes Volume and $E^*$ the set of analogs of $E$’s members. 

\medskip

However, consider rewriting Area by the following (mathematically equivalent) expression:
\[
\text{\emph{Area}}^{**}(x ,y)= x^{2-1}y^{2-1}-\sin^{2-2}(x\cdot y).
\]
Since we know that the square maximizes \emph{Area}$^{**}$, one might reason that the cube maximizes:
\[
\text{\emph{Volume}}^{**}(x, y, z)=x^{3-1}y^{3-1}z^{3-1}-\sin^{3-2}(x\cdot y \cdot z).
\]
The analogical inference in this case seems extremely weak. Yet, since $Q(H \vert E)$ is high (where H is re-written in terms of \emph{Area}$^{**}$) by Bartha’s symmetry argument one ought to assign a ‘close enough’ value to the bet $Q (H^{**} \vert E^*)$, where $H^{**}$ is the claim that the cube maximizes \emph{Volume}$^{**}$ and $E^*$ is the same set of analogues of $E$’s members as before. Yet, this is highly implausible: it should be at least permissible for a rational agent to assign negligible value to the latter bet.\footnote{Note that the similarity between \emph{Area}$^{**}$ and \emph{Volume}$^{**}$ counts as ‘admissible’ on Bartha’s (2009:170--1) account since it (unlike the case that Bartha 2009:169 presents) preserves the symmetry of the argument’s source proof.}

\medskip

To summarize, we find that Bartha’s account is too weak in one sense and too strong in another. It is too weak because it neglects cases in which an argument from analogy supports a mathematical conjecture that is already assigned non-negligible probability (or betting quotient), as in the case of Deligne’s proof indirectly supporting the Riemann hypothesis. It is also too strong because, as the result of its commitment to a symmetry-based justification, it ends up giving confirmatory power (in the non-incremental ‘plausibility’ sense that Bartha identifies) to analogical arguments in mathematics that should not be understood as having such power at all. The appeal to symmetry thus flattens all analogical arguments in mathematics, preventing us from recognizing those that are evidentially relevant to a conjecture, and even investigating why they are so. For these reasons, we think that a more refined framework is necessary. The next section assesses the prospects for a Bayesian approach to confirmation by mathematical analogy.

\section{Bayesianism in Mathematics}
\label{Sec3}

In his classic discussion on analogy and induction in mathematics, Polya (1954b) indicates BC  as a potential candidate for capturing the role of analogical arguments in making plausible new mathematical conjectures. His proposal (which is discussed approvingly in Corfield 2003:83) is that confirmation in mathematics comes in different degrees and that the degree of support offered by an analogical argument is proportional to the “hope” (Polya 1954b:27) that a “common ground” (27) exists between the analogy’s source and target. In other words, the idea is that an analogical argument is stronger when we may regard the similarities figuring in the argument as pointing to some deep but yet unknown connection between the two mathematical domains. Although Polya does not develop this idea in full probabilistic detail, his proposal appears to respond to at least one of the requirements that we found inadequately addressed by Bartha’s account: viz., the idea that different analogical arguments provide different degrees of confirmation to mathematical conjectures. As Corfield (2003) perceptively notes, the appeal to one’s expectation for a common ground promises to offer suitable ground for such a distinction. 

\medskip

Notwithstanding Polya’s and Corfield’s optimism, the prospects for applying the formal tools of BC to the case of confirmation by mathematical analogy are far from straightforward. Two varieties of obstacles need to be overcome, neither of which is trivial. A general problem is logical omniscience (mentioned above). Following Garber (1983), we can distinguish two readings of BC’s aims. On the “thought police model” (1983:101), BC’s aim is to patrol our reasoning as we collect more evidence. On the “learning machine model”, instead, BC’s aim is to describe the intellectual life of an ideal learner. In both cases, the coherence requirement – viz., the idea that degrees of belief should obey the axioms of probability – forces one to have credence equal to $1$ in any logical truth and credence of $0$ in any contradiction. Moreover, the same axioms require that the probability of a deductive conclusion should not be smaller than the probability of the conjunction of the premises. It follows that anyone who (no matter how mathematically skilled) is not logically omniscient is thereby irrational from BC’s standpoint. 

\medskip

There is no agreed-upon way of solving the problem of logical omniscience. One well-known attempt is by Garber (1983). The guiding idea is to preserve logical omniscience for all propositions true in virtue of their truth-functional form (tautologies), while relaxing omniscience for all propositions that follow validly from others (not true in virtue of their truth-functional form). This idea is implemented by one’s treating  “$A \vdash B$” in a given propositional language L (over which the probability function is defined) as atomic propositions. Thus, we replace all instances of “$A \vdash B$” in $L$ severally with some ‘dummy’ expression $\tau$ and apply the probability function over the newly obtained language $L^*$. In this way, the logical implication of $A$ to $B$ becomes ‘invisible’ to the Bayesian apparatus, permitting an agent to be rationally uncertain as to whether $\tau$ holds. As Eells (1990) notes, however, the solution is incomplete, since an agent is still required to assign credence $1$ to all tautologies. This is already an implausible demand on rationality since some tautologies are computationally very complex. 

\medskip

A more selective approach, suggested by Gaifman (2004), consists in applying the probability function to only a subset of (the formulas of) the language of interest. The intuition is to limit the application of probability to that part of the language that is easily grasped by a realistic agent (e.g., because it contains relatively short deductions) and to add specific axioms (regarding implication and the notion of “local provability”, which is not transitive), so that logical omniscience can hold solely within the selected subset of the given language. While this proposal overcomes the incompleteness affecting Garber’s proposal, it is not exempt from defects of its own. For instance, as Easwaran (2009) notes, the proposal seems unable to regard as rationally permissible a mathematician’s uncertainty about statements of which she has proposed a proof - one that she perfectly understands in all of its steps. Such uncertainty is bound to be deemed as irrational on Gaifman’s proposal, although it sometimes appears to be rational. 

\medskip

n light of the problems above, we may take a more pragmatic approach and regard logical omniscience as simply one more of BC’s yet unresolved issues - alongside, e.g., the well-known “problem of old evidence” (Glymour 1980).\footnote{Another option is to insist, with Franklin (2013; 2020), that confirmatory relations exist and are knowable independently of entailments. On Franklin’s view, this case is exactly analogous to that of an agent who assigned a credence of $0.5$ to a coin landing heads despite knowing that the world in which she lives is fully deterministic. Neither the proposal nor the comparison are, however, fully convincing: it seems that it is at least sometimes irrational to hold an intermediate degree of belief in a proposition when one is in a position to prove (or disprove) it.} Indeed, this is precisely the approach that we are going to assume below. Dialectically, asking for such a waiver is not question-begging. After all, we have at least some partial proposals as to how to relax logical omniscience. Moreover, better attempts are not unlikely to become available with more time and inventiveness. In the meantime, we know that BC has been fruitfully applied to inductive reasoning in the empirical sciences. Accordingly, we may regard as our central concern that of understanding what a Bayesian approach can tell us about confirmation in mathematics. This will determine whether it is worth spending time carefully reflecting on more foundational issues - in particular, whether the hard problem of relaxing logical omniscience deserves our further intellectual efforts.

\medskip

Even if a waiver for logical omniscience is granted, specific issues about the application of BC to analogical reasoning remain open. It is, after all, one thing to claim that an analogical argument may contribute to increasing one’s trust in a conjecture from BC’s standpoint; it is another to demonstrate this claim formally by means of BC’s probabilistic apparatus. Presumably, any proposed representation would need to display some sort of diagnostic value: it would need to demonstrate that it yields a verdict of confirmation (or lack thereof) solely in virtue of the fact that a strong (respectively, weak) analogical argument is at work. Without some such diagnostic capacity, we will not say that BC properly captures analogical reasoning. It is not obvious that BC can achieve that. Furthermore, even if BC can properly capture some examples of analogical inference in mathematics, it is not obvious that it can capture all of them. Neither Polya (1954b) nor Corfield (2003) provide sufficient details about the intended application of BC to judge if their proposals are sufficiently illuminating and comprehensive in the above respects. 

\medskip

To summarize, BC’s reliance on degrees of confirmation promises great descriptive fit with mathematical practice. However, the prospects for a Bayesian analysis of analogical inferences in mathematics are all but obvious. For one thing, there is the problem of logical omniscience. But even setting that general problem aside, there remains the specific (and highly non-trivial) problem of providing the details of the promised Bayesian account of confirmation by analogy. In the next section, we are going to develop a proposal about representing inductive support from analogy within BC. Even though we believe that such a proposal represents a useful addition to understanding a variety of cases in which analogies play a role in mathematical research, our ultimate recommendation will be to embrace a middle ground in between Bartha’s position, on one hand, and Polya and Corfield’s, on the other. As section five will contend, a hybrid framework - one which recognizes both incremental and non-incremental aspects of confirmation by analogy in mathematics - holds the greatest promise of satisfying all desiderata.

\section{Incremental Confirmation by Analogy}
\label{Sec4}

Let a reasonable notion of probability as applied to mathematical conjectures be given. In other words, let’s assume that the problem of logical omniscience is solved by one of the responses discussed above (it does not matter which one). Here is one thing that we can hope to model by means of a Bayesian framework: roughly, the \emph{incremental} form of confirmation that attaches to a hypothesis about the target as a consequence of the discovery that some analogous result holds in the source. The case from which we shall start is that of the confirmation of Riemann’s hypothesis due to Deligne’s proof of Weil’s analogue conjecture. As it turns out, a Bayesian representation of such a case is possible that is at least in the spirit of Polya’s (1954b) suggestion concerning the ‘hope for a common ground’. Let’s see what the actual formalization looks like.

\subsection{The Riemann Hypothesis}
\label{sec41}

Some historical background about the case-study will be useful. Riemann’s (1859) hypothesis emerged in the context of studying prime numbers distribution. Using a procedure called \emph{analytic continuation}, he generalized to the complex plan a function introduced by Euler in connection with his solution to the so-called ‘Basel problem’. Thus, the \emph{Riemann zeta function} was defined:
\begin{equation}
\label{e4}
\zeta(s)=\sum_{n=1}^{\infty}\frac{1}{n^s},
\end{equation}      		   
where the sum is over the natural numbers and s is a complex variable different from one. The main result that Riemann established was a relation between the zeros of his zeta functions and the distribution of prime numbers. Some of these zeros are called \emph{trivial} in that their existence  easy to prove (they are exactly the negative even numbers). The question is about the \emph{non-trivial} ones. It is known that they belong to the open strip of the complex plane defined by the set of complex numbers with real part between $0$ and $1$. The Riemann hypothesis is that all non-trivial zeros lie on the critical line defined by the complex number having their real part equal to $\frac{1}{2}$.

\medskip

While resisting several attempts at a proof, another area of mathematics unexpectedly came to the rescue of the Riemann hypothesis: algebraic geometry. Its principal object, the ‘algebraic variety’, is defined as the set of solutions to a system of polynomials. The coefficients of such polynomials belong to so-called ‘finite fields’, a set whose operations are analogous to those of the real numbers, but containing only a finite number of elements. Furthering Artin’s (1924) studies into this relatively recent field, André Weil (1949) demonstrated that it is possible to construct a zeta function involving the number of points over the finite extension of the original field. One of Weil’s four conjectures about this \emph{algebraic} zeta function mimics the Riemann hypothesis in considering the distribution of the zeros of the integral polynomials linked to the rational form of the algebraic zeta function. Building upon Weil’s work, Deligne (1974) eventually proved the analogue of Riemann hypothesis for algebraic varieties over finite fields.

\medskip

Soon after his algebraic formulation, Weil began to ponder about the bearing of his conjecture to the original number-theoretical question as formulated by Riemann:
\begin{quote}
The Riemann hypothesis [...] appears to-day in a new light, which shows it to be closely connected with the conjecture of Artin on the L-functions, thus making these two problems two aspects of the same arithmetico-algebraic question. (1950:297)
\end{quote}
The same sentiment of a connection between the two conjectures is shared by those several experts in the field who regard Deligne’s proof as “the best reason to believe that [the Riemann hypothesis] is true” (Deninger 1994:493). Their reasoning is based on an appreciation of the analogy between the domains that Riemann’s and Weil’s works respectively explored. This “sense of the relatedness of mathematical facts” (Corfield 2003:121) underwrites their allegations to the effect that the Riemann hypothesis has been confirmed by Deligne’s proof.

\medskip

Our claim is that the experts’ reasoning can be represented formally as an instance of the phenomenon of transitivity of confirmation. As it is well-known, the Bayesian notion of confirmation is not generally transitive: sometimes $A$ confirms $B$ and $B$ confirms $C$, but it is not the case that $A$ confirms $C$. However, there exist formal proposals in the philosophical literature (e.g., Hesse 1970; Roche and Shogenji 2013) about when transitivity of confirmation is guaranteed to occur. We propose that confirmation of the Riemann hypothesis can be understood, from a Bayesian perspective, as an instance of the case in which the weakest conditions for transitive confirmation hold. To illustrate this Bayesian analysis, we define:
\begin{itemize}
    \item[($R$)] The Riemann hypothesis is true.
\end{itemize}

We then make precise the evidence that, by the lights of several experts, supports R by analogy:
\begin{itemize}
    \item[($W$)] There is a proof of Weil’s conjecture.
\end{itemize}
At this point, it can be shown that Bayesian confirmation of $R$ by $W$ occurs just in case non-extremal credence (i.e. neither zero nor one) is given to the following bridge claim:\footnote{A credence of one to $G$ means that there is no uncertainty about the analogy and hence R can be deduced directly from the known facts about the target. Cf. also Weil (1979:251): “we have the theory [..., so] gone is the analogy”.} 
\begin{itemize}
    \item[($G$)] The distribution of the solutions of Weil’s zeta function is robust (i.e., does not change) in the passage from the algebraic-geometric question to the analogous number-theoretic question,
\end{itemize}
and the following probabilistic conditions between $R$, $W$, and $G$ are jointly satisfied:

\begin{itemize}
    \item[a)] $\P(R \vert G) > \P(R)$;
    \item[b)] $\P(W \vert G) > \P(W \vert \lnot G)$;
    \item[c)] $	\P(R \vert G \wedge W) \ge \P(R \vert G)$;
    \item[d)] $\P(R \vert \lnot G \wedge W) \ge \P(R \vert \lnot G)$.
\end{itemize}

We note that all of a)-d) are plausibly satisfied. Specifically, conditions a) and b) plausibly express one’s disposition (which is based on the “sense of relatedness of the mathematical facts” discussed above) to increase trust in the Riemann hypothesis as the result of the discovery of an analogous result in algebraic geometry. This is expressed in terms of one’s seriously entertaining a ‘bridge’ hypothesis $G$, defined in such a way as to connect two mathematical domains from the standpoint of confirmation. Conditions c) and d) are weak additional assumptions, serving to ensure that $G$’s confirmation appropriately transfers to $R$. They require merely that $W$ (that there is a proof of Weil’s conjecture) does not disconfirm $R$ conditional on knowing $G$ (the presence of a common ground) or, alternatively, on knowing $\lnot G$ (the absence of a connection). Both assumptions are plausible. In particular, d) is justified on grounds that $W$ is arguably irrelevant to $R$ if $\lnot G$ is true; consequently, $\P(R \vert \lnot G \wedge W) = \P(R \vert \lnot G)$.\footnote{Note that c) is trivial if we take $G$ to entail $R$. Insofar as we are trying to be neutral on the response to the problem of logical omniscience, we include c) so as to avoid assuming that the agent can ‘see’ the implication of $R$ by $G$.} A theorem of transitive confirmation by Roche and Shogenyi (2013; see Appendix I) shows that a)-d) together entail $\P(R \vert W) > \P(R)$, meaning that there is confirmation of $R$ by $W$ in the standard Bayesian sense. 

\medskip

In summary, the above model shows that, provided it makes sense to assign non-extremal probabilities to mathematical hypotheses such as $R$, $W$ and $G$, it is possible to spell out a precise Bayesian account of confirmation of the Riemann hypothesis from the proof of Weil’s analogue conjecture.\footnote{To be clear, it is not our claim here that any rational agent ought to accept conditions a)-d). We recognize that there may be mathematicians whose background knowledge justifies different assignments of probabilities from those stated. For such agents, confirmation of $R$ by $W$ would simply not occur. The point of the above Bayesian model is to show that, under a reasonable distribution of probabilities - one that plausibly expresses an expert’s recognition of the ‘relatedness’ of number theory and algebraic geometry - there can be confirmation in BC’s sense.} Such a full-fledged Bayesian formalization - a true novelty for the existing literature - vindicates Polya’s insight that the degree of confirmation depends upon the ‘hope’ that a common ground exists between the mathematical domains being compared. The proposal above shows that Polya’s insight is correct at least insofar as confirmation depends upon one’s disposition to accept a non-extremal credence to a hypothesis ($G$) acting as a confirmatory bridge between source and target. As a way of validating our proposed Bayesian approach to incremental confirmation of mathematical conjectures, let’s now turn to another case-study.  

\subsection{Taylor expansion of functions}
\label{sec42}

Here is another case of confirmation that can be treated in a Bayesian fashion. It is the case in which the discovery of a striking \emph{similarity} between two apparently unrelated mathematical domains confirms a hypothesis about the target. (This is different from the case in \ref{sec41}, where the discovery of a result in a source supports a conjecture in an analogous target; see fn. 13). Here is an elementary illustration (cf. Lange 2021 for an alternative take on the example). By calculating the Taylor series around $0$ of the function $\frac{1}{1-x^2}$, we note that it is similar to that of $\frac{1}{1+x^2}$: 
\begin{align*}
\frac{1}{1-x^2}&=1+x^2+x^4+x^6+\cdots;\\
\frac{1}{1+x^2}&=1-x^2+x^4-x^6+\cdots.
\end{align*}

An additional aspect of similarity is that the respective series have similar convergence behavior: both converge when $|x|<1$ and diverge when $|x|>1$. This is so despite the fact that the corresponding functions have different shapes in the real plane (in particular, $\frac{1}{1-x^2}$ is not defined for $x=1$). From these similarities, a mathematician might infer that $\frac{1}{1+x^2}$ is related to $\frac{1}{1-x^2}$ by some connection that is simply not ‘visible’ in the real plane. This (factually correct) reasoning is based on the idea that it would be something of a coincidence if their respective series just happened to have similar expressions and convergence behavior, but were otherwise unrelated.

\medskip

Using BC, we can represent the effect on a rational agent’s credences of the discovery of the similarities between the two functions. To obtain an adequate Bayesian model, the first step is to clearly define the evidence introduced by the analogical inference. This is arguably that:
\begin{itemize}
    \item[($O$)] The Taylor series of the source and target functions display similar expressions and convergence behavior.
\end{itemize}
We then precisely define the conclusion about the target that the argument aims to support:
\begin{itemize}
    \item[($T$)] The function $\frac{1}{1+x^2}$ is related to $\frac{1}{1-x^2}$ by a connection not visible in the real plane.
\end{itemize}
Finally, we define an appropriate ‘bridge’ between source and target (more on this below):
\begin{itemize}
    \item[($J$)] Some deeper mathematical fact holds such that the two functions are connected.
\end{itemize}
By Roche and Shogenji’s (2013) transitivity theorem, when $O$, $T$ and $J$ are assigned non-extremal credence, then O confirms T just in case the following conditions obtain:
\begin{itemize}
    \item[e)] $\P(T \vert J) > \P(T)$;
    \item[f)] $\P(O \vert J) > \P(O \lvert \lnot J)$;
    \item[g)] $\P(T \vert J \wedge O) \ge \P(T \vert J)$;
    \item[h)] $\P(T \vert \lnot J \wedge O)  \ge \P(T \vert \lnot J)$. 
\end{itemize}
The conditions plausibly describe the epistemic situation of the expert mathematician who comes to suspect a connection between the two functions ($T$) as the result of some newly observed similarity between them ($O$). In particular, e) and f) plausibly express one’s defeasible recognition that the similarities in the respective series would be somewhat unlikely if there were no ‘common ground’ (existing in some deeper mathematical domain) relating the two functions ($J$). Conditions g) and h) ensure that confirmation of the bridge $J$ (the common ground) by the similarities $O$ properly translates into $T$’s confirmation. Both assumptions are plausible. In particular, h) is arguably justified on grounds that $O$ is evidentially irrelevant to $T$ if $J$ is false; accordingly, $\P(T \vert \lnot J \wedge O) = \P(T \vert \lnot J)$. It follows from acceptance of e)-h) that $\P(T \vert O) > \P(T)$. In informal terms, the model shows that there is confirmation of $T$ by the observed similarities.

\medskip

One point worth noting here is that, while formally the conditions are the same as in case \ref{sec41}, the content of the evidence and the bridge hypothesis has changed in the new example. This is not an \emph{ad hoc} move. It is justified by the fact that the cases of confirmation by analogy in \ref{sec41} and \ref{sec42} are importantly different. In the former case, we have two mathematical domains that we already suspect to be connected in some way; the similarities between the two domains constitute the background information in virtue of which, from the discovery of a result in the source, one confirms the corresponding conjecture in the target. In the latter case, we have instead the discovery of a striking similarity between two apparently unrelated domains (the two functions); here the similarities constitute the evidence in virtue of which the hypothesis of a common ground is confirmed (insofar as it makes the similarities more likely than they would be if there were no common ground).\footnote{This distinction between two types of analogical inferences in mathematics is already noted in Polya (1971:42), who names the former (exemplified by \ref{sec41}) an inference from “\emph{methods}” to “\emph{results}” and the latter an inference (exemplified by \ref{sec42}) from “\emph{results}” to “\emph{methods}”. Note that both types of analogical inferences in mathematics are fueled by one’s ‘hope for a common ground’ - although, as stressed above, they do so in slightly different ways.   } Both cases exemplify equally valid, though distinct notions of confirmation by analogy in mathematics, which the account above correctly treats differently. 

\subsection{The Euler Characteristic}
\label{sec43}
The distinction above can be brought to bear on an important test for the proposed Bayesian analysis: the case-study of the Euler characteristics examined in Example 2.2. As mentioned earlier, the case is interesting because two different analogical arguments can be used to support the conclusion that a regularity exists uniting the elements of a convex polyhedron (expressed by $V-E+F-S=1$), analogous to that which holds for regular polygons in plane geometry (i.e., $V=E$). The two arguments arguably differ with respect to the support that they provide to the conclusion - the Euler characteristic. If a Bayesian representation is to be a useful way of representing confirmation by analogy in the realm of mathematics, then, it should be able to tell apart the two analogical arguments with respect to their respective capacity to provide inductive support. (We noted earlier that the above is a \emph{desideratum} that Bartha’s 2009 rival account fails to satisfy). Let’s take a look as to how well the present Bayesian analysis fares in this regard.

\medskip

Based on the distinction drawn in \ref{sec42}, we should note immediately that the two analogical arguments differ not just in strength but in kind. The Cauchy-inspired argument resembles case \ref{sec41} in that it aims to support its conclusion by using as background information the geometrical similarities that link all convex polyhedra together (in virtue of which they belong to the same geometrical genus); from these assumed (mostly tacit) similarities, and the fact that smashed cubes, tetrahedra, and dodecahedra satisfy the Euler characteristic, one defeasibly infers that the latter holds true for any convex polyhedron whatsoever. Conversely, the algebraic, Polya-inspired argument is similar in kind to the inference in case \ref{sec42}: from the similarity between the expression $V-E+F-S=1$ (known to hold for cubes, tetrahedra and polyhedra) and the corresponding two-dimensional equation $V-E+F=1$, one defeasibly infers the existence of a deeper mathematical fact that is capable of explaining the observed similarity.

\medskip

Based on this preliminary observation, we fix the conclusion that both arguments arrive at:
\begin{itemize}
    \item[($C$)] All convex polyhedra satisfy equation \eqref{e1} - the so-called ‘Euler characteristic’.
\end{itemize}
Let us then reconstruct the two arguments in accordance with the framework developed in \ref{sec41} and \ref{sec42}. For the Cauchy-inspired reasoning, the new evidence introduced by the argument is that:
\begin{itemize}
    \item[($E$)] Smashed cube, tetrahedra, and dodecahedra satisfy the property whereby the relation $V-E+F$  is constant even after removal of any inner triangular surface.
\end{itemize}
The bridge hypothesis is formulated in accordance with a recipe indicated for case \ref{sec41}:
\begin{itemize}
    \item[($B$)] The result of the smashing, triangularization and removal operations on convex polyhedra is robust with respect to the type of the polyhedron.
\end{itemize}

As in case \ref{sec41}, by ‘robust’ here we mean that the result remains invariant under changes in the property specified in the bridge hypothesis - in this case, a change in the form of the polyhedron.

\medskip

For the Polya-inspired argument, both the evidence and the bridge hypothesis need to be defined differently. The new evidence that the argument makes salient is the observation that:
\begin{itemize}
    \item[($E^*$)] Cubes, tetrahedra, and dodecahedra satisfy an alternate sum regularity \eqref{e3} algebraically similar to that which unites the elements of regular polygons \eqref{e2}.
\end{itemize}
The bridge hypothesis is formulated in accordance with the recipe for case \ref{sec42}, as follows:
\begin{itemize}
    \item[($B^*$)] Some deeper mathematical fact holds such that \eqref{e2} is algebraically similar to \eqref{e3}.
\end{itemize}	
The ‘deeper mathematical fact’ in this case may be a general theorem (in some yet unknown mathematical theory) having as a consequence that the alternating sum of the number of elements of an n-dimensional figure that satisfies some precise constraints is always equal to $1$.

\medskip

We are interested in comparing the analogical component of the confirmation of $E$ and $E^*$ to $C$. On our account, the conditions for confirmation by analogy are, respectively, i)-l) and m)-p):

\begin{multicols}{2}
\begin{itemize}
    \item[i)] $\P(C \vert B) > \P(C)$;
    \item[j)] $\P(E \vert B) > \P(E \vert \lnot B)$;
    \item[k)] $\P(C \vert B \wedge E) \ge \P(C \vert B)$;
    \item[l)] $\P(C \vert \lnot B \wedge E)  \P(C \vert \lnot B)$.
\end{itemize}

\columnbreak

\begin{itemize}
    \item[m)] $\P(C\vert B^*) > \P(C)$;
    \item[n)] $\P(E^* \lvert B^*) > \P(E^* \lvert \lnot B^*)$;
    \item[o)] $	\P(C \vert B^* \wedge E^*) \ge \P(C \vert B^*)$;
    \item[p)]	$\P(C \vert \lnot B^* \wedge E^*)  \ge \P(C \vert \lnot B^*)$. 
\end{itemize}
\end{multicols}

As in the previous cases \ref{sec41} and \ref{sec42}, we assume that a reasonable distribution of the agent’s probabilities, whereby each of $E$, $E^*$, $B$, $B^*$ and C are assigned non-extremal prior credence. Let’s now consider the content of the probabilistic conditions for each of the two columns above.

\medskip

The conditions in the first column, i)-l), are all plausibly satisfied. Specifically, i) and j) arguably express one’s recognition that, as the result of the geometrical similarities that link all convex polyhedra together, similar manipulations (smashing, triangularization, etc.) are likely to lead to the same result: in compressing a polyhedron, some initially unconnected vertices will be found on the same surface and will thus be potential vertices for (roughly) triangular figures; the removal of any such triangular surfaces is likely to keep the number of elements constant. Conditions k) and l) are weak additional assumptions, both of which are plausibly satisfied. In particular, l) can be justified on grounds that E is evidentially irrelevant to $C$ if $\lnot B$ is known. In more informal terms, if we know that the result of similar manipulations is not robust to the smashing and triangularization operations across different polyhedra, then E arguably ceases to be evidence for $C$. Altogether, acceptance of i)-l) entails the correct verdict that $\P(C \vert E) > \P(C)$.

\medskip

Conversely, it is not obvious that m)-p) are all satisfied. Condition n) is especially shaky. Considering the algebraic similarity between \eqref{e2} and \eqref{e3}, the expert mathematician (which we are assuming not to know if $C$ holds) may well regard it as little or no evidence for $C$. After all, the similarity with \eqref{e3} may appear to be the accidental fallout of an algebraic manipulation of \eqref{e1} into \eqref{e2}, one which has no independent basis in geometry. Accordingly, a mathematician may reject $E^*$’s capacity to confirm (and not merely suggest) that some deep mathematical fact underlies the algebraic similarity between \eqref{e2} and \eqref{e3}, and therefore that all convex polyhedra (and not merely the observed ones) obey the alternating sum regularity that \eqref{e3} expresses. Hence,  $E^*$ will fail to confirm $C$ in anything but a weak enumerative sense. That is to say, even though $E^*$ may rationally increase the probability of the generalization $C$, it would only do so to a much slighter degree than the confirmation achieved by means of the geometrical reasoning. This is exactly as we should expect given the perceived difference in strength between the arguments.\footnote{To clarify, we are not claiming that the Polya-inspired argument fails to provide a significant degree of confirmation (for a contextually relevant threshold of significance) because it is an argument from results to methods (see fn. 13); rather, it fails because some salient aspects of the relevant contextual information fail to back it up.}

\bigskip

To wrap up, we have now spelled out a precise proposal about capturing confirmation by mathematical analogy in Bayesian terms. We have offered recipes for two kinds of incremental confirmation by analogy, respectively exemplified in \ref{sec41} and \ref{sec42}, showing that both can be viewed as instances of transitive confirmation.\footnote{Hesse’s (1970) discussion on transitivity of confirmation and its relation to analogy is an early inspiration for the present proposal. See also Nappo (2021) for a recent revival of Hesse’s idea in the domain of the empirical sciences.} As we have seen in the (somewhat artificial, but still instructive) case of the Euler characteristic, the proposed framework possesses a fair degree of diagnostic capacity. Specifically, when fed with the relevant background knowledge (understood in terms of reasonable credence assignments to conjectures), it can help pinpoint the factors that make certain inferences by analogy stronger than others. By identifying the precise features of the background knowledge that determine whether, and to what extent, confirmation by analogy occurs, the Bayesian proposal on offer can thus be an important addition to our understanding of analogy in mathematics. The next section will address the question of what the account above omits with regards to the practice of reasoning from analogy in mathematics.

\section{Non-Incremental Confirmation}
\label{Sec5}
The previous section has argued that, assuming that the problem of logical omniscience can be solved, we can represent many important instances of confirmation by mathematical analogy by means of standard BC. Even though the sample of case studies is necessarily limited, the illustrations provided in the previous section at least show that the notion of confirmation by analogy does not violate the tenets of BC and can thus be understood as a fully legitimate part of the inductive methodology of mathematics from a Bayesian standpoint. Moreover, as we stressed earlier, the recipes that we have provided for two types of analogical inferences in mathematics, exemplified by cases \ref{sec41} and \ref{sec42}, reflect independently plausible ideas about when an analogical argument in mathematics is inductively strong. This fact makes us reasonably confident that the proposed framework is appropriately generalizable beyond the illustrations considered above.

\medskip

In this section, our aim is to show that there is one notion of confirmation by analogy that the Bayesian account above is, by its nature, not prepared to capture - what we call the ‘non-incremental’ notion of confirmation by analogy. What we will be primarily concerned with clarifying is how this notion interacts with the incremental form just discussed. For an illustration, let’s consider the geometrical case of the analogical inference from area to volume examined in Example 2.3. From the fact that, of all rectangles, the square maximizes \emph{Area}, we are led to expect that, of all rectangular boxes, the cube maximizes \emph{Volume}. The argument can be represented formally by means of the recipe indicated above. As in case \eqref{sec41}, we have the discovery of a result in a source (viz., that the square maximizes \emph{Area}) which bears on a conjecture about solid geometry (viz., that the cube maximizes \emph{Volume}).\footnote{To avoid complications, we assume that the evidence collected in the source is new (i.e., not already part of the background knowledge). For an insightful discussion on the problem of ‘old evidence’ for BC, see Sprenger (2015).} On our proposal, confirmation in BC’s sense occurs if an agent’s credence function countenances a ‘bridge’ hypothesis to the effect that the result about squares is robust to the passage to solid geometry.

\medskip

One aspect of the Bayesian formalization (which we leave as an exercise for the reader, being a simple extension of the recipe in \ref{sec41}) must be emphasized. It assumes that a Bayesian agent would find it reasonable to assign non-negligible probability to the ‘bridge’ hypothesis that links squares to cubes from the standpoint of confirmation. Although such an assumption seems eminently reasonable in light of our grasp of geometrical relations, it must be stressed that nothing in the Bayesian framework dictates this choice. As far as the axioms of probability are concerned, there is no difference between epistemic agents who take it that cubes are the three-dimensional analogue of squares and agents who take it that (say) spheres are the three-dimensional analogue of squares. The latter agents may find it natural to think that the sphere’s surface area is the analogue of a square’s area but not that a cube’s volume is. Presumably, such deviant agents may fail to even entertain the bridge hypothesis that we find natural in this context, viz., that which makes cubes the three-dimensional analogue of squares.\footnote{The problem echoes (and is closely related to) Goodman’s (1954) ‘new riddle of induction’; cf. Franklin (2013).}

\medskip

We are, of course, considering an extreme example for illustration. But the issue is much closer to home than one might suppose. Some mathematicians, for instance, categorically reject the analogy between Weil’s theorems for algebraic varieties and Riemann’s hypothesis in number theory (cf. Corfield 2003:98). Presumably, such experts - a minority among those who can be regarded as fully grasping the mathematical problem - assign little or no credence to the bridge hypothesis that links the proof of Weil’s conjecture for algebraic varieties to the Riemann hypothesis. Accordingly, they may regard Deligne’s proof of Weil’s conjecture as little or no evidence for the truth of the Riemann hypothesis.\footnote{Cf. also Corfield (2003) on the “monstrous moonshine” connection between the J-function and the ‘monster group’ (a problem in modular representation theory): “if one had asked a mathematician how likely she thought it that there be some [...] connection between them [...], the answer would presumably have been ‘vanishingly small’” (125).} Here we have a realistic case in which a difference in the bridge hypotheses that an agent considers overturns the verdict of confirmation.  

\medskip

Reflecting on these cases is useful because they suggest the existence of two distinct notions of confirmation by analogy. On the one hand, there is the ‘incremental’ sense in which the discovery of some new fact about a source, or some additional similarity between a source and a target, may justify \emph{additional} credence to a conjecture about the target. That is the notion that the Bayesian account of the previous section aims to explicate. On the other hand, there is the ‘non-incremental’ sense in which a hypothesis is judged to be plausible in light of an analogy between a source and a target mathematical domain.\footnote{The distinction between two notions of confirmation parallels Carnap’s (1950), who distinguished roughly between the ‘incremental’ sense in which an empirical hypothesis is confirmed by an additional piece of evidence and the ‘absolute’ sense in which a hypothesis may be regarded as highly plausible in light of one’s total evidence.} For instance, as the natural analogue of squares in solid geometry, we find it plausible that cubes will satisfy the analogue of various geometrical properties that hold for squares. This is not the sense of confirmation in which we can say that a piece of evidence confirms (to some degree) a hypothesis. What we have instead is an all-thing-considered \emph{judgment} to the effect that a given hypothesis (what we have been calling a ‘bridge’ hypothesis) deserves non-negligible credence. In this sense, it is non-incremental.\footnote{It must be stressed that the distinction that we are advancing is not psychological. With regards to the question of how mathematicians actually process analogical reasoning, we defer to experts in cognitive science.}  

\medskip

It would be a mistake to suppose that, since the judgments of analogy in question do not contribute to the incremental confirmation of a mathematical conjecture, they cannot make a difference (in some broad sense) to a mathematician’s epistemic life. In many circumstances, mathematicians may fail to note the analogy between two domains of mathematical interest; or, while vaguely recognizing some element of resonance, they may fail to pinpoint exactly what the two domains have in common (cf. Polya 1954b:111). Under those circumstances, a clearer recognition of the connection between two domains may rationally lead a mathematician to consider ‘bridge’ hypotheses that she had not even considered before. The result, from a Bayesian standpoint, would be a revision of the agent’s credence function to include a new hypothesis. The difference with the incremental form of confirmation by mathematical analogy is that such change in credence would not be dictated by the consideration of any new evidence. Rather, we can say that such change derives from looking at one’s old evidence in a new way: as possibly indicating a previously unnoticed connection between distinct mathematical domains.  

\medskip

A question that may be pressed at this point is whether the recognition of a non-incremental notion of confirmation by analogy is consistent with the adoption of a Bayesian epistemology. On this issue, we have little doubt: there is no tension. Since BC is silent about which combinations of priors and likelihoods that preserve coherence with the axioms of probability should be adopted by a rational agent, it is perfectly consistent for a Bayesian to recognize that judgments of analogy may rationally inform an epistemic agent’s priors and, at the same time, to hold that the only rational way for an agent to update her credences on the basis on new evidence is by applying the standard diachronic rule of conditionalization. Consistency with BC remains even if one accepts the plausible doctrine that there is often a fact of the matter as to which judgments of analogy in particular contexts of research actually justify consideration of which bridge hypotheses.\footnote{ So-called ‘Dutch-book’ objections (cf. Bartha 2009:299) do not arise against this proposal since we are denying that there is any set of rules of reasoning the application of which necessarily yields the ‘correct’ judgment of analogy in a given context. See also Bartha (2019:$\S$5.1) on the position that he labels “liberal Bayesianism”.} Insofar as judgments of analogy of the kind we are pointing to pertain to the sphere of prior credences, our entitlement to them systematically escapes BC’s field of concern.

\medskip

Of course, our claim that, even before new similarities and dissimilarities are considered, judgments of analogy may play a role in determining the ‘reasonable’ probabilities that a mathematician assumes in the course of research shows that there is an intrinsic limit to Bayesian reconstructions of analogical reasoning in mathematics. Even though we can appeal to BC’s probabilistic apparatus to precisely describe the mechanisms whereby a mathematician updates her credences on the basis of new similarities and dissimilarities, questions with regards to which ‘bridge hypotheses’ should be given consideration, and precisely to what extent, are not answered from within BC. Because of this, we may wonder whether any rival formal framework may be able to do better than a Bayesian approach. If that were the case, the recognition of a non-incremental form of confirmation by analogy would be an argument against BC’s adoption.

\medskip

However, we are highly skeptical that there can be such a rival framework. Even though the judgments of analogy in question are not in principle unanalyzable, we find it hard to imagine where to even start to provide a formal account of them, i.e., a theory that successfully reduces the norms governing the appropriateness of certain judgments of analogy in mathematical practice to the application of some general rules of reasoning. To give a simple illustration of the difficulties, consider the plausibility judgments involved in another geometrical example (cf. Bartha 2009:110). It is a theorem of plane geometry that the three medians of a triangle intersect in a point, called a ‘centroid’. By a ‘median’, we mean the segment that unites the midpoint of a side to the opposite vertex. One conjecture in solid geometry that the theorem induces is the following (Fig. 2): the four medians of a tetrahedron similarly intersect in a point. In this case, by a tetrahedron’s ‘median’ we mean the segment that unites a face’s centroid to the opposite vertex. 

\bigskip

\bigskip

\begin{figure}[ht]
\centering
\begin{minipage}[c]{.49\textwidth}
\centering
\begin{tikzpicture}[scale=1.25,>=latex]
\draw[ultra thick] (0,0)--(4,0)--(1,3)--(0,0);
\draw[red] (2,0)--(1,3);
\draw[red] (0.5,1.5)--(4,0);
\draw[red] (2.5,1.5)--(0,0);
\draw[ultra thick] (0,0)--(4,0);
  \fill (2,0) circle[radius=3pt];
  \fill (0.5,1.5) circle[radius=3pt];
  \fill (2,0) circle[radius=3pt];
  \fill (2.5,1.5) circle[radius=3pt];
  \fill[red] (2.5/3*2,1.5/3*2) circle[radius=3pt];
\end{tikzpicture}
\end{minipage}
\begin{minipage}[c]{.49\textwidth}
\centering
\begin{tikzpicture}[scale=1.25,>=latex]
\draw[ultra thick] (4,0)--(1.5,3.5);
\draw[red] (1.5 , 1.25) --(1,2);
\draw[red] (2.167,1.833)--(-1.5,1);
\draw[red] (1.167,1)--(1.5,3.5);
\draw[red] (0.333,2.167)--(4,0);
  \fill (1.167,1) circle[radius=3pt];
  \fill (0.333,2.167) circle[radius=3pt];
  \fill (2.167,1.833) circle[radius=3pt];
  \fill (1.5 , 1.25) circle[radius=3pt];
  \fill[red] (1.25,1.66) circle[radius=3pt];
  \draw[ultra thick] (-1.5,1)--(4,0)--(1,2)--(-1.5,1);
\draw[ultra thick] (-1.5,1)--(1.5,3.5)--(1,2);
\end{tikzpicture}
\end{minipage}
\caption*{Figure 2}
\end{figure}

\newpage

Although the conjecture about tetrahedra is highly plausible, the strength of the analogical inference depends on geometrical intuitions that are difficult to articulate. Possibly the easiest way to express them consists in borrowing the physical notion of a figure’s ‘center of mass’. The median of a triangle can be regarded as the point where the ideal mass of the triangle concentrates. Similarly, the median of a tetrahedron is where the mass of the tetrahedron (assuming it had any mass and that it were uniformly distributed across its volume). Given that the ‘median’ in both the plane and the solid case is what unites the center of mass of a figure or face to the opposite vertex, and given that we know triangles and tetrahedra to be analogous to one another in other ways, it is plausible that results about medians in triangles will have an analogue for tetrahedra. Insofar as making this judgment requires a complex combination of visualization and knowledge of the relevant geometrical properties of triangles and tetrahedra, it is very hard to imagine what combination of algorithmic procedures of reasoning would be able to capture the judgments of the trained mathematician even in such an elementary example.

\medskip

The example above illustrates how, even though a mathematician’s training and practice often predisposes her to see certain connections among mathematical domains rather than others, and to be accordingly induced to perform certain inferences rather than others, exactly which connections and inferences is not determined by a rule. For one, these judgments often rely on mental visualization, which is hard to explicate. Moreover, mathematical creativity often consists precisely in the capacity of seeing new connections between mathematical domains when the latter were initially thought to be separate. In Polanyi’s (1958) famous terms, we may consider this a form of ‘tacit knowledge’ - a hard-to-articulate inferential know-how that comes with mathematical training. The capacity to see such connections and to separate them into deeper and more superficial ones is thus something more akin to an \emph{aesthetic} capacity - as when developing a capacity to recognize more and less beautiful pieces of art and to make one own’s beautiful pieces- rather than to an algorithmic procedure of reasoning. Because of this, we are highly skeptical that any formal account will fare significantly better than BC.

\medskip

In summary, our claim is that we can capture one notion of confirmation by analogy in Bayesian terms but not another. Even though the judgments of analogy that underlie the non-incremental notion are not universally agreed upon, they are often widely shared - to the point that it becomes plausible to assume that the acquisition of the proper training for mathematical research is measured at least in part by one’s disposition to express certain judgments of analogy in the course of research rather than others. Yet, besides showing that the judgments of analogy of a ‘reasonable’ mathematician are encapsulated in assignments of non-zero credence to certain bridge hypotheses rather than others, it seems not plausible to expect that a formal framework such as BC will be able to say something more informative on the subject. At the same time, as we stressed, the recognition of a non-incremental notion of confirmation by analogy does not pose trouble for Bayesians insofar as analogy’s role in informing priors is consistent with the adoption of BC as a doctrine about credence updating. 

\medskip

A comparison with Newton’s theory of gravity may be helpful at this point. It is well-known that Newton formulated his law of universal gravitation without advancing hypotheses as to the ultimate cause of the force of gravity. Our hybrid framework for analogical reasoning in mathematics does something similar. We have shown that (assuming that logical omniscience has been appropriately relaxed) it is possible to represent in Bayesian terms the incremental form of confirmation that attaches to a mathematical conjecture as the result of the discovery of an additional fact about the source, or an additional similarity between source and target domains, given a reasonable distribution of prior and likelihoods. But as to the role that analogy plays in supporting that initial distribution, \emph{hypotheses non fingimus}. Although the judgments of analogy in question may be sufficiently robust and stable across contexts to allow for plausible generalizations to be made about them, we are unlikely to ever reduce a trained mathematician’s capacity to form plausible judgments of analogy to the application of some formal rules. It follows that the hybrid account gives us all we could hope for in terms of precise formalization.

\section{Conclusion}
\label{Sec6}

In this paper, we have provided a framework that systematically addresses the inferential mechanism whereby evidence in a familiar mathematical domain supports conjectures in an analogous target. For this purpose, we have identified two notions of confirmation by analogy that are relevant to the realm of mathematics – an \emph{incremental} and a \emph{non-incremental} notion - and discussed how each of their roles is to be understood from the standpoint of BC. While acknowledging logical omniscience as a yet unresolved issue for BC, one of our major achievements consists in having provided the details of a Bayesian account of the incremental notion of confirmation by analogy, one which (given an answer to the logical omniscience problem) responds to questions about precisely when an analogy may contribute to increasing a mathematician’s credence in the truth of a conjecture, as well as precisely what sorts of facts about the contextual information are relevant to showing that there is (or is not) confirmation.   

\medskip

We have defended our account as a superior alternative to those proposed by, respectively, Bartha (2009) and Corfield (2003). With the help of accessible case-studies from mathematical practice, we have contended that Bartha’s (2009) proposal is too strong in one sense and too weak in another. It is too strong because it makes too many analogies capable of providing conjectures with a form of ‘prima facie’ plausibility; it is also too weak because it fails to account for the different degrees of support that analogical arguments can provide to conjectures that already possess some ‘non-negligible’ probability (in Bartha’s sense). With regards to Corfield’s (2003) Bayesian approach, we have mainly pointed out its failure to specify the details of a Bayesian account of confirmation by analogy. As we have seen in the previous section, once the details are fully explicit it becomes clear that there are two distinct forms of confirmation by analogy to account for and that, while a Bayesian account of the incremental form can be provided, a formal approach is unlikely to capture the the non-incremental form of confirmation.

\medskip

Of course, much work remains to be done to understand the role of analogy in mathematics. Among other things, the relation between analogy in mathematics and in the empirical sciences requires further investigation. More attention also needs to be paid to the normative problem of articulating a justification for the use of analogical reasoning in the mathematical domain. Although our proposal directly responds to neither of these further issues, we believe it points towards new and fruitful answers - ones that we will systematically develop in future work.

\section*{Appendix I}

Roche and Shogenji’s (2013) transitivity theorem states the following:
\begin{theorem} 
\label{T1}
Given any three propositions with non-extremal credences $X$, $Y$ and $Z$, if:
\begin{itemize}
    \item[(i)] $\P(Z \vert Y) > \P(Z)$
    \item[(ii)] $\P(X \vert Y) > P (X \vert \lnot Y)$
    \item[(iii)] $\P(Z \vert X \& Y) \ge \P(Z \vert Y)$
    \item[(iv)] $\P(Z \vert X \& \lnot Y) \ge \P(Z \vert \lnot Y)$
\item[then:]
    \item[(v)] $ \P(Z \vert X) <\P(Z)$.
\end{itemize}

\end{theorem}

\begin{corollary}
In the limiting case in which $Y$ entails $Z$, (ii) and (iv) guarantee (v).
\end{corollary}

\noindent
Please refer to Roche and Shogenji (2013) for a proof of Theorem \ref{T1}.

\end{document}